\documentclass[10pt,twoside]{article}

\newcommand{\sub}{\supseteq}

\def\ifif {if and only if\ \ }
\def\sub {\subseteq}

\usepackage{amsmath}
\usepackage{amsfonts}
\usepackage{amssymb}
\usepackage{makeidx}
\usepackage{amsthm}
\newcommand{\mathscr}{\mathbf}

\begin{document}
\title{Intersections of essential minimal prime ideals}
%\dedication{\langle Professor Azarpanah \rangle}
%\author{ F. Azarpanah \\Department of Mathematics, Chamran University, Ahvaz , Iran}
%\thanks{2000 {\it Mathematics Subject Classification}: 54C40.}
\author{A. Taherifar\\ Department of Mathematics, Yasouj University, Yasouj,
Iran\\}
\date{ataherifar@mail.yu.ac.ir}
 \maketitle
%\maketitle
%--------------------------------------------------------
\theoremstyle{definition}\newtheorem{thm}{Theorem}[section]
\theoremstyle{definition}\newtheorem{cor}[thm]{Corollary}
\theoremstyle{definition}\newtheorem{lem}[thm]{Lemma}
\theoremstyle{definition}\newtheorem{prop}[thm]{Proposition}
\theoremstyle{definition}\newtheorem{defn}[thm]{Definition}
\theoremstyle{definition}\newtheorem{Rem}[thm]{Remark}
\theoremstyle{remark}\newtheorem{rem}[thm]{Remark}
\theoremstyle{definition}\newtheorem{exam}[thm]{Example}
\theoremstyle{definition}\newtheorem{qu}[thm]{question}
\newtheorem{ques}[thm]{Question}
%---------------------------------------------------------------------

{\noindent\bf Abstract.} Let $\mathcal{Z(R)}$ be the set of zero
divisor elements of a commutative ring $R$ with identity and
$\mathcal{M}$ be the space of minimal prime ideals of $R$ with
Zariski topology. An ideal $I$ of $R$ is called strongly dense
ideal or briefly $sd$-ideal if $I\subseteq \mathcal{Z(R)}$ and is
contained in no minimal prime ideal. We denote by
$R_{K}(\mathcal{M})$, the set of all $a\in R$ for which
$\overline{D(a)}=\overline{\mathcal{M}\setminus V(a)}$ is compact.
We show that $R$ has property $(A)$ and $\mathcal{M}$ is compact
\ifif $R$ has no $sd$-ideal. It is proved that
$R_{K}(\mathcal{M})$ is an essential ideal (resp., $sd$-ideal)
\ifif $\mathcal{M}$ is an almost locally compact (resp.,
$\mathcal{M}$ is a locally compact non-compact) space. The
intersection of  essential minimal prime  ideals of a reduced ring
$R$ need not be an essential ideal. We find an equivalent
condition for which any (resp., any countable) intersection of
essential minimal prime ideals of a reduced ring $R$ is an
essential ideal. Also it is proved that the intersection of
essential minimal prime ideals of $C(X)$ is equal to the socle of
C(X) (i.e., $C_{F}(X)=O^{\beta X\setminus I(X)}$). Finally, we
show that a topological space $X$ is pseudo-discrete \ifif
$I(X)=X_{L}$ and $C_{K}(X)$ is a pure ideal.

%------------------------------------------------------------------------------------%
\vspace{0.3 cm}\noindent{\it Keywords:} Essential ideals,
$sd$-ideal, almost locally compact space, nowhere dense, Zariski
topology.

\vspace{0.3 cm}\noindent{\it 2010 Mathematical Subject
Classification: $13A15$, $54C40$}
%------------------------------------------------------------------------------------%
\section{Introduction}
In this paper, $R$ is assumed  to be a commutative ring with
identity, and $\mathcal{M}$ is the space of minimal prime ideals
of $R$. For $a\in R$ let $V(a)=\{P\in\cal{M}$: a$\in P\}$. It is
easy to see that for any $R$, the set $\{D(a)=\mathcal{M}\setminus
V(a): a\in R \}$ forms a basis of open sets on $\cal{M}$. This
topology is called the Zariski topology. For $A\sub R$, we use
$V(A)$ to denote the set of all $P\in\cal{M}$ where $A\sub P$ (see
\cite{8.}). For a subset $H$ of $\cal{M}$ we denote by
$\overline{H}$ the closure points of $H$ in $\mathcal{M}$. The
intersection of all minimal prime ideals containing $a$ is denoted
by $P_{a}$ . An ideal $I$ of $R$ is called a $z^{0}$-ideal, if
$P_{b}\subseteq P_{a}$ and $a\in I$, then $b\in I$ (See \cite{2.}
and \cite{5.}). For any subset $S$ of a ring $R$, ${\rm
ann}$$(S)=\{a\in R: aS=0\}$.
\\\indent We denote by $C(X)$ the ring of real-valued,
continuous functions on a completely regular Hausdorff space X,
$\beta X$ is the Stone-Cech compactification of $X$ and for any
$p\in \beta X$, $O^{p}$ (resp., $M^{p}$) is the set by all $f\in
C(X)$ for which $p\in int_{\beta X}cl_{\beta X} Z(f)$ (resp.,
$p\in cl_{\beta X} Z(f)$). Also, for $A\sub\beta X$, $O^{A}$ is
the intersection of all $O^{p}$ where $p\in A$. It is well known
that $O^{p}$ is the intersection of all minimal prime ideals
contained in $M^{p}$. We denote the socle of $C(X)$ by $C_{F}(X)$;
it is characterized in \cite{11.} as the set of all functions
which vanish everywhere except on a finite number of points of X.
The known ideal $C_{K}(X)$ in $C(X)$, is the set of functions with
compact support, and the generalization of this ideal is defined
in \cite{14.}. The reader is referred to \cite{7.} for undefined
terms and notations.\\\indent A non-zero ideal in a commutative
ring is said to be essential if it intersects every non-zero ideal
non-trivially, and the intersection of all essential ideals, or
the sum of all minimal ideals, is called the socle (see
\cite{12.}).
\\\indent We denote by $R_{K}(\mathcal{M})$, all  $a\in R$ for
which $\overline{D(a)}$ is compact as a subspace of $\mathcal{M}$.
In section $2$, by algebraic properties of the ideal
$R_{K}(\mathcal{M})$, we find topological properties of the space
of minimal prime ideals of $R$, $\mathcal{M}$, say locally
compactness and almost locally compactness. Also we call $I$  a
strongly dense ideal or briefly $sd$-ideal if $I\subseteq
\mathcal{Z(R)}$ and is contained in no minimal prime ideal. We
characterize commutative reduced rings $R$ which have no
$sd$-ideals. It is proved that $R_{K}(\mathcal{M})$ is contained
in the intersection of all strongly dense $fd$-ideals (i.e., if
${\rm ann}(F)\subseteq{\rm ann}(a)$, for some finite subset $F$ of
$I$ and $a\in R$, then $a\in I$).\\\indent In section $3$, we show
that the intersection of essential minimal prime ideals in any
ring need not be an essential ideal. In a reduced ring $R$, we
prove that every intersection of essential minimal prime ideals is
an essential ideal if and only if the set of isolated points of
$\mathcal{M}$ is dense in $\mathcal{M}$. Also it is proved that
every countable intersection of essential minimal prime ideals of
a reduced ring $R$ is an essential ideal if and only if every
first category subset of $\mathcal{M}$ is nowhere dense in
$\mathcal{M}$. We characterize the intersection of essential
minimal prime ideals in $C(X)$, i.e., the intersection of
essential minimal prime ideals of $C(X)$ is equal to the ideal
$C_{F }(X)$ (i,e., the socle of $C(X)$). Finally, we prove that
the intersection of essential minimal prime ideals of $C(X)$ is
equal to the ideal $C_{K}(X)$ \ifif $I(X)=X_{L}=\bigcup_{f\in
C_{K}(X)}\overline{X\setminus Z(f)}$, i.e., $I(X)=X_{L}$ and
$C_{K}(X)$ is a pure ideal. By this result and Theorem $4.5$ in
\cite{3.}, we see that $X$ is a pseudo-discrete space \ifif
$I(X)=X_{L}$ and $C_{K}(X)$ is a pure ideal.

%------------------------------------------------------------------------------------%

\section{ $R_{K}(\mathcal{M})$ and Strongly dense ideals}

In this section we introduce the ideal $R_{K}(\mathcal{M})$ and
the class of strongly dense ideals as a subclass of dense ideals.
We show that a reduced ring $R$ has no $sd$-ideal \ifif $T(R)$
(i.e., the total quotient ring of $R$) is a von Neumann regular
ring. By this, we have $C(X)$ has no $sd$-ideal \ifif $X$ is a
cozero-complemented space. It is proved that $R_{K}(\mathcal{M})$
is an essential ideal (resp., $sd$-ideal) \ifif $\mathcal{M}$ is
an almost locally compact space (resp., locally compact
non-compact space).
\begin{defn}
Let  $R$ be a commutative ring with identity and $D(a)$ be the set
of all prime ideals which do not contain $a$. We define the family
$R_{K}(\mathcal{M})$ to be the set of all $a\in R$ for which
$\overline{D(a)}$ is compact (as a subspace of $\mathcal{M}$).
\end{defn}
%--------------------------------------------------------
\begin{exam}
If $\cal{M}$ is a discrete space, then $R_{K}(\mathcal{M})=\{a\in
R: D(a)\ is finite\}$. For example, let $R$ be the weak (discrete)
direct sum of countably many copies of the integers. $R$ may be
regarded as the ring of all sequences of integers that are
ultimately zero. Then $\cal{M}$ is a countable discrete space (see
\cite[2.11]{8.}).
\end{exam}
%-----------------------------------------------------------
\begin{lem}\label{1.1}
\begin{itemize}
\item[\rm(i)] $R_{K}(\mathcal{M})$ is a $z^{0}$-ideal of $R$.
\item[\rm(ii)] $R_{K}(\mathcal{M})=R$ \ifif $\mathcal{M}$ is
compact. \item[\rm(iii)] $R_{K}(\mathcal{M})=0$ \ifif
$\mathcal{M}$ is nowhere compact (i.e., the interior of every
compact set is empty).
\end{itemize}
\end{lem}
\textbf{Proof.} (i) For $a,b\in R_{K}(\mathcal{M})$, we have
$\overline{D(a+b)}\subseteq \overline{D(a)}\cup\overline{ D(b)}$,
and if $a\in R$, $b\in R_{K}(\mathcal{M})$, then
$\overline{D(ab)}\subseteq \overline{D(b)}$. Therefore
$R_{K}(\mathcal{M})$ is an ideal of $R$. Now let $P_{b}\subseteq
P_{a}$ and $a\in R_{K}(\mathcal{M})$. Then $V(a)\subseteq V(b)$,
hence $\overline{D(b)}\subseteq \overline{D(a)}$ so
$\overline{D(b)}$ is compact, i.e., $b\in R_{K}(\mathcal{M})$.

(ii) By definition, it is obvious.

(iii) $\Rightarrow$ Let $K$ be compact subset of $\mathcal{M}$ and
$P\in int(K)$ Then there is a non-zero element $f\in R$ such that
$P\in D(f)\subseteq int(K)$, so $f\in R_{K}(\mathcal{M})=0$, i.e.,
$D(f)=\phi$, which is a contradiction.\\\indent $\Leftarrow$
Suppose that  $f\in R_{K}(\mathcal{M})$. Then $D(f)$ is contained
in the interior of $\overline{D(f)}$ so $D(f)=\phi$, i.e., $f=0$.
$\hfill\square$
%-----------------------------------------------------------------------
\begin{defn}
An ideal $I$ of $R$ is called a strongly dense ideal or briefly a
$sd$-ideal if $I\subseteq \mathcal{Z(R)}$ and is contained in no
minimal prime ideal ($V(I)=\phi$).
\end{defn}
%---------------------------------------
\begin{exam}
\begin{itemize}
\item[\rm(i)] Every  prime ideal of a ring $R$ which is not a
minimal prime and contained in $\mathcal{Z(R)}$ is a
$sd$-ideal.\item[\rm(ii)] For a set $X$, let $R=\Bbb R^{X}$ (i.e.,
the ring of real valued functions). Then we can see that any
element of $R$ is an unit or a zero-divisor. So any ideal $I$ of
$R$ for which $V(I)=\phi$ is a $sd$-ideal.\item[\rm(iii)] Let $p,
q$ be two non-isolated points in an almost $P$-space $X$ (i.e.,
every zero-set has nonempty interior (see \cite{L} and \cite{V})).
Then the ideal $I=M_{p}\cap M_{q}$ is a $sd$-ideal.
\end{itemize}
\end{exam}
%-----------------------------------------------------------
Recall that an ideal $I$ of $R$ is a dense ideal if
$\rm{ann}$$(I)=0$. We observe that in any commutative reduced ring
$R$ the ideal $I\oplus\rm{ann}$$(I)$ is an essential ideal. Hence
an ideal $I$ in a reduced ring $R$ is an essential ideal \ifif it
is a dense ideal \cite{12.}. \\\indent In the following, we see
that a non-minimal prime ideal need not be a $sd$-ideal.
\begin{Rem}\label{1.2}
 Every $sd$-ideal in a reduced ring $R$ is a dense ideal (essential ideal), but there is a dense ideal which is not $sd$-ideal.
 To see this, let $I$ be a $sd$-ideal and $g\in \rm{ann}$$(I)$. Then $gf=0$ for each $f\in I$ therefore we have, $\bigcap_{f\in I} V(gf)= \mathcal{M}$,
 hence $V(g)\cup(\bigcap_{f\in I} V(f))=\mathcal{M}$, i.e., $V(g)=\mathcal{M}$, thus $g=0$. Now let $x$ be a non-isolated point in compact space $X$.
 Then by \cite[Remark 3.2]{4.}, the ideal $O_{x}$ is an essential ideal of $C(X)$ which is not a minimal prime ideal. By \cite[Theorem, 3.1]{4.}, ${\rm
 ann}$$(O_{x})=0$, i,e., $O_{x}$ is dense ideal.
 However, this ideal is not a $sd$-ideal.
 Because there is a minimal prime ideal in $C(X)$ contains $O_{x}$, i.e., $V(O_{x})\neq\phi$.
\end{Rem}
%------------------------------------------------------------------
We denote by $I_{z}$, the intersection of all $z$-ideals that
contain $I$. An ideal $I$ of $R$ is called a $rez$-ideal if there
is an ideal $J$ for which $I\not\subseteq J$ and $I_{z}\cap
J\subseteq I$. For more see \cite{2.} and \cite{5.}.
\begin{prop} Every ideal $I$ in a reduced ring $R$,  is a $rez$-ideal or  a dense ideal.
\end{prop}
\textbf{Proof.} Let $I$  be a non-$rez$-ideal in $R$. By \cite[
Corollary 2.8]{2.}, $\rm{ann}$$(I)=0$, so $I$ is a dense ideal.
$\hfill\square$
%----------------------------------------------------------------
\begin{lem}\label{1.3}
 Let $R$ be a reduced ring.
 \begin{itemize}
\item[\rm(i)] $\bigcap_{i=1}^{n}V(f_{i})=\phi$ \ifif
$\bigcap_{i=1}^{n}\rm{ann}$$(f_{i})=0$. \item[\rm(ii)] If $F$ is a
finite subset of $R$, then $V(F)=\mathcal{M}\setminus
V(\rm{ann}$$(F))$. \item[\rm(iii)] If $I\subseteq \mathcal{Z(R)}$
is a finitely generated ideal, then $I$ is an $sd$-ideal \ifif $I$
is a dense ideal. \item[\rm(iv)]If $R$ has  finitely many minimal
prime ideals, then $R$ has no $sd$-ideal.
\end{itemize}
\end{lem}
\textbf{Proof.} Trivial. $\hfill\square$
%---------------------------------------------------------------------
\\ \indent
Recall that a ring $R$ has property $(A)$ (resp., property
$(a.c.)$), if for every finitely generated ideal $I\subseteq
\mathcal{Z(R)}$, $\rm{ann}$$(I)\neq 0$ (resp., for any finitely
generated ideal $I$ of $R$ there is $c\in R$ such that
$\rm{ann}$$(I)=\rm{ann}$$ (c)$), see \cite{8.} and \cite{H}.
\\\indent In the following theorem we characterize a class of reduced rings
which has no $sd$-ideal.
\begin{thm}\label{1.5} Let $R$ be a reduced ring with total quotient $T(R)$. The following conditions are equivalent.
\begin{itemize}
\item[\rm(i)] $T(R)$ is a von Neumann regular ring.
\item[\rm(ii)]$R$ satisfies property $(A)$ and $\mathcal{M}$ is
compact.\item[\rm(iii)] $R$ has no $sd$-ideal. \item[\rm(iv)] $R$
satisfies property $(a.c)$ and $\mathcal{M}$ is compact.
\end{itemize}
\end{thm}
\textbf{Proof.} For (i)$\Leftrightarrow$(ii)$\Leftrightarrow$(iv),
see \cite[Theorem 4.5]{10.}.\\\indent (ii)$\Rightarrow$(iii) Let
$I$ be a $sd$-ideal. Then $I\sub \cal{Z}(R)$ and $\bigcap_{f\in
I}V(f)=\phi$. Hence $\mathcal{M}=\bigcup_{f\in I}D(f)$.
 Compactness of $\mathcal{M}$ implies that there are $f_{1},...,f_{n}\in I$ such that $\bigcap_{i=1}^{n}V(f_{i})=\phi$. By Lemma
 \ref{1.3}, we have
 $\bigcap_{i=1}^{n}\rm{ann}$$(f_{i})=\rm{ann}$$(F)=0$, where $F=\{f_{1},...,f_{n}\}$. This is a contradiction, for $R$ has property $(A)$ and $F\subseteq
 \mathcal{Z(R)}$.\\\indent (iii)$\Rightarrow$(ii)
 Suppose that $I\subseteq \mathcal{Z(R)}$ is a finitely generated ideal and $\rm{ann}$$(I)=0$. Then by Lemma \ref{1.3},
 $I$ is a $sd$-ideal, this  contradicts  the hypothesis. Thus $R$ has property $(A)$. Now, let $\mathcal{M}=\bigcup_{f\in S}D(f)$ where $S$ is a
 proper subset of $R$.
 If $(S)\subseteq \mathcal{Z(R)}$,
 then $\bigcap_{f\in (S)}V(f)=\phi$ implies that the ideal generated by $S$ is a $sd$-ideal, a contradiction.
  Hence $(S)\not\subseteq \mathcal{Z(R)}$, so there is $H=\{f_{1},...,f_{n}\}\subseteq S$ such that $\rm{ann}$$(H)=0$. By Lemma \ref{1.3},
  $V(H)=\bigcap_{i=1}^{n}V(f_{i})=\mathcal{M}\setminus V(\rm{ann}$$(H))=\emptyset$, thus $\mathcal{M}=\bigcup_{i=1}^{n}D(f_{i})$, i.e., $\cal{M}$ is compact. $\hfill\square$
  %-----------------------------------------------------------------------------------------------
\\\indent
Henriksen and  Woods have introduced cozero complemented spaces
(i.e., if, for every cozero-set $V$ of $X$, there is a disjoint
cozero-set $V^{\prime}$ of $X$ such that $V\cup V^{\prime}$ is a
dense subset of $X$) (see \cite{9.}). In \cite{8.}, they  have
proved that the space of minimal prime ideals of $C(X)$ is compact
\ifif $X$ is a cozero complemented space. Now by Theorem
\ref{1.5}, and the fact that $C(X)$ satisfies property $(a.c)$ we
have the following corollary.
%-----------------------------------------------
\begin{cor}
$C(X)$ has no $sd$-ideal \ifif $X$ is a cozero complemented space.
\end{cor}
%---------------------------------------------------------------------------------
Recall that a ring $R$ has property $(c.a.c)$, if for any
countably generated ideal $I$ of $R$, there exists $c\in R$ such
that $\rm{ann}$$(I)=\rm{ann}$$(c)$, see \cite{8.}. If $R$ is a
ring with property $(c.a.c)$, then by \cite[Theorem, 4.9]{8.}
$\mathcal{M}$ is countably compact. But in a ring with property
$(A)$ this need not be true.
\begin{prop}
 Let $R$ be a reduced ring. $R$ satisfies property $(A)$ and $\mathcal{M}$ is countably compact
\ifif $R$ has no countably generated $sd$-ideal.
\end{prop}
\textbf{Proof.} The proof is similar  to that of the proof of
Theorem \ref{1.5} step  by step.$\hfill\square$
%--------------------------------------------------------------------
\\ \indent
Recall that an ideal $I$ of $R$ is called an $fd$-ideal, if for
each finite subset $F$ of $ I$ and $x\in R$, ${\rm
ann}$$(F)\subseteq{\rm ann}$$(x)$ implies that $x\in I$. For more
details see \cite{13.}.
%--------------------------------------------------------------------
\begin{prop}\label{1.6}
\begin{itemize}
\item[\rm(i)] $R_{K}(\mathcal{M})$ is contained in the
intersection of all strongly dense  $fd$-ideals in $R$.
\item[\rm(ii)] $R_{K}(\mathcal{M})$ is a $sd$-ideal \ifif
$\mathcal{M}$ is a locally compact non-compact space.
\end{itemize}
\end{prop}
\textbf{Proof.} (i) Let $I$ be a strongly dense $fd$-ideal, $f\in
R_{K}(\mathcal{M})$ and $P\in \overline{D(f)}$. Then there is
$g\in I$ such that $P\in D(g)$, so $\overline{D(f)}\subseteq
\bigcup_{g\in I}D(g)$. On the other hand $\overline{D(f)}$ is
compact so there are $g_{1},...,g_{n}\in I$ such that
$\overline{D(f)}\subseteq \bigcup_{i=1}^{n}D(g_{i})$. Hence
$V(f)\supseteq \bigcap_{i=1}^{n}V(g_{i})=V(\{g_{1},...,g_{n}\})$.
This implies that $\rm{ann}$$(F)\sub \rm{ann}$$(f)$, where
$F=\{g_{1},...,g_{n}\}\subseteq I$. But  $I$ is a $fd$-ideal so
$f\in I$. \\\indent (ii) Let $R_{K}(\mathcal{M})$ is a $sd$-ideal
and  $P\in \mathcal{M}$. By definition, there is $f\in
R_{K}(\mathcal{M})$ such that $P\in D(f)\subseteq \overline{D(f)}$
so $P$ has a compact neighborhood, i.e., $\mathcal{M}$ is a
locally compact space. On the other hand, $\mathcal{M}$ is not
compact, since if $\mathcal{M}$ is compact, then by Lemma
\ref{1.1}, $R_{K}(\mathcal{M})=R$, which is a
contradiction.\\\indent Conversely, first, we see that
$R_{K}(\mathcal{M})\subseteq Z(R)$. Because, if $f\in
R_{K}(\mathcal{M})$ and $\rm{ann}$$(f)=0$, then $D(f)=\mathcal{M}$
is compact, which is a contradiction, by hypothesis. Now for every
$P\in \mathcal{M}$ there is a compact neighborhood $K$ of $P$ in
$\mathcal{M}$. So there is $f\in R$ such that $P\in D(f)\subseteq
intK\subseteq K$, i.e., $f\in R_{K}(\mathcal{M})$, hence
$R_{K}(\mathcal{M})$ is an $sd$-ideal. $\hfill\square$
%---------------------------------------------------------------------------
\\\indent
A Hausdorff space $X$ is said to be an almost locally compact
space if every non-empty open set of $X$ contains a non-empty open
set with compact closure (see \cite{3.}). The next result is a
topological characterization of $R_{K}(\mathcal{M})$ as an
essential ideal.
%---------------------------------------------------------------------
\begin{thm}\label{1.7}
Let $R$ be a reduced ring. Then $R_{K}(\mathcal{M})$ is an
essential ideal \ifif $\mathcal{M}$ is an almost locally compact
space.
\end{thm}
\textbf{Proof.} $\Rightarrow$ Let $R_{K}(\mathcal{M})$ be an
essential ideal and $U$ be an open subset of $\mathcal{M}$. Then
there exists a non-zero element $f\in R$ such that $D(f)\subseteq
U$. It is enough to prove that $D(f)$ contains $D(g)$ for some
$g\in R_{K}(\mathcal{M})$. If $D(f)\cap D(g)=\phi$ for each $g\in
R_{K}(\mathcal{M})$, then $D(fg)=\phi$, so $fg=0$, i.e.,
$R_{K}(\mathcal{M})\cap (f)=0$, this is a contradiction by
essentiality of $R_{K}(\mathcal{M})$. Hence there is $g\in
R_{K}(\mathcal{M})$ such that $D(fg)=D(f)\cap D(g)\neq\phi$, but
$D(fg)\subseteq D(f)$, i.e., $U$ contains an open subset for which
the closure is compact.\\\indent $\Leftarrow$ Let $f$ be a
non-zero element in $R$. It is enough to prove that
$R_{K}(\mathcal{M})\cap (f)\neq\phi.$ $D(f)\neq \phi$ is an open
subset in $X$. By hypothesis there is an open subset $V\subseteq
D(f)$ such that $\overline{V}$ is compact, so there is a non-zero
element $g\in R$ such that $D(g)\subseteq V\subseteq D(f)$, i.e.,
$g\in R_{K}(\mathcal{M})$. Now $D(fg)=D(f)\cap D(g)=D(g)\neq\phi$,
hence $fg\neq 0$ is an element of $R_{K}(\mathcal{M})\cap (f)$,
i.e., $R_{K}(\mathcal{M})$ is an essential ideal. $\hfill\square$
%-------------------------------------------------------

\section{ Intersections of essential minimal prime ideals}

The intersection of  essential minimal prime  ideals of a reduced
ring $R$ need not be an essential ideal. Even a countable
intersection of essential minimal prime ideals need not be an
essential ideal. For example, the ideal $O_{r}$ for any rational
$0\leq r \leq 1$ is an essential ideal in $C(\Bbb R)$, which is
the intersection of minimal prime ideals. Now for any $0\leq r
\leq 1$ let $P_{r}$ be a minimal prime ideal that contains
$O_{r}$. Then any $P_{r}$ is an essential ideal, but by
\cite[Theorem, 3.1]{3.}, the ideal $I=\bigcap P_{r}$ is not an
essential ideal, for $\bigcap Z[I] = [0, 1]$ and $int[0, 1]\neq
\phi$. In this section we give a topological characterization of
the intersection of essential minimal prime ideals of a reduced
ring $R$ (resp., $C(X)$) which is an essential ideal.\\\indent For
an open subset $A$ of $\cal{M}$, suppose that $O_{A}:=\{a\in R:
A\sub V(a)\}$. Since for any $a, b\in R$, $V(a)\cap V(b)\sub
V(a-b)$ and for each $r\in R$, $a\in O_{A}$, we have $V(a)\sub
V(ra)$,  thus $O_{A}$ is an ideal of $R$. It is easy to see that
$O_{A}=\bigcap_{P\in A}P$ and $V(O_{A})=\overline{A}$, where
$\overline{A}$ is the cluster points of $A$ in $\cal{M}$.\\\indent
 We need the following lemmas which are easy to prove.
%------------------------------------------------------------------------------------%
\begin{lem}\label{1.8}
Let $R$ be a reduced ring. An ideal $I$ of $R$ is an essential
ideal \ifif $intV(I)=\phi$.
\end{lem}
%------------------------------------------------------------------
\begin{lem}\label{l} The intersection of all essential minimal prime ideals in a reduced ring $R$ is equal to the ideal $O_{(\mathcal{M}\setminus I(\mathcal{M}))}$,
where $I(\mathcal{M})$ is the set of isolated points of
$\mathcal{M}$.
\end{lem}
%--------------------------------------------------------------
In \cite{3.}, Corollary 2.3 and Theorem 2.4, Azarpanah showed that
every intersection (resp., countable intersection) of essential
ideals of $C(X)$ is essential \ifif the set of isolated points of
$X$ is dense in $X$ (resp., every first category subset of $X$ is
nowhere dense in $X$). Now, we generalize these results for the
essentiality of the intersection of essential minimal prime ideals
in a reduced ring.
%----------------------------------------------------------------
\begin{prop}\label{1.9} Let $R$ be a reduced ring.
Every intersection of essential minimal prime ideals is an
essential ideal \ifif the set of isolated points of $\mathcal{M}$
is dense in $\mathcal{M}$.
\end{prop}
\textbf{Proof.} Assume that every intersection of essential
minimal prime ideals is an essential ideal. Then Lemma \ref{l}
implies that $O_{(\mathcal{M}\setminus I(\mathcal{M}))}$ is an
essential ideal. By Lemma \ref{1.8},
$intV(O_{(\mathcal{M}\setminus I(\mathcal{M})})=\phi$. On the
other hand, we have
\begin{center}
$V(O_{\mathcal{M}\setminus
I(\mathcal{M})})=\overline{\mathcal{M}\setminus
I(\mathcal{M})}=(\mathcal{M}\setminus I(\mathcal{M}))$.
\end{center}
Therefore $int(\mathcal{M}\setminus I(\mathcal{M}))=
int(V(O_{\mathcal{M}\setminus I(\mathcal{M})}))=\phi$. This shows
that $\overline{I(\mathcal{M})}=\mathcal{M}$.\\\indent Conversely,
by hypothesis, $int(V(O_{\mathcal{M}\setminus
I(\mathcal{M})}))=int(\mathcal{M}\setminus I(\mathcal{M}))=\phi$.
So by Lemma \ref{1.8}, $O_{(\mathcal{M}\setminus I(\mathcal{M}))}$
is an essential ideal. Since $O_{(\mathcal{M}\setminus
I(\mathcal{M}))}$ is contained in  every intersection of essential
minimal prime ideals, so  every intersection of essential minimal
prime ideals is an essential ideal. $\hfill\square$
%---------------------------------------------------------------------------
\begin{thm}
Let $R$ be a reduced ring.
 Every countable intersection of essential minimal prime ideals of $R$ is an essential ideal
 \ifif every first category subset of $\mathcal{M}$ is nowhere dense in $\mathcal{M}$.
\end{thm}
\textbf{Proof.} $\Rightarrow$
 Let $(F_{n})$ be a sequence of nowhere dense subsets of $\mathcal{M}$. Then  by Lemma \ref{1.8}, for each $n\in \Bbb
 N$, the ideal
 $O_{F_{n}}=\bigcap_{P\in F_{n}}P$,  is an essential ideal. By hypothesis,
 $E=\bigcap_{n=1}^{\infty}O_{F_{n}}=O_{(\bigcup_{i=1}^{\infty}F_{n})}$ is an essential ideal. On the other hand
 $V(E)=\overline{(\bigcup_{i=1}^{\infty}F_{n})}$. So we must have $int(\overline{\bigcup_{n=1}^{\infty}{F_{n}}})=\phi$, i.e., $\bigcup_{i=1}^{\infty}F_{n}$ is nowhere dense.
\\ \indent
 $\Leftarrow$ Let $(I_{n})$ be a sequence of essential minimal prime ideals in $R$. Letting $\{I_{n}\}=F_{n}$,
 then $int(F_{n})=intV(I_{n})=\phi$, i.e., each $F_{n}$ is a nowhere dense subset of $\mathcal{M}$. $O_{F_{n}}\subseteq I_{n}$ implies that
 $O_{A}\subseteq \bigcap_{n=1}^{\infty}I_{n}$, where $A=\bigcup_{n=1}^{\infty}F_{n}$. Now we have $V(O_{A})=\overline{A}$, and since $A$ is a first category subset,
 then $int(\overline{A})=\phi$, i.e., $O_{A}$ is an essential ideal. Thus $\bigcap_{n=1}^{\infty}I_{n}$
 is also an essential ideal. $\hfill\square$
%-------------------------------------------------------------
\\\indent
The following lemma is a characterization of  the intersection of
all essential minimal prime ideals of $C(X)$.
%--------------------------------------------------------------------
\begin{lem}\label{1.10}
The intersection of all essential minimal prime ideals of $C(X)$
is equal to $O^{\beta X\setminus I(X)}$,  where $I(X)$ is the set
of isolated  points of topological space $X$.
\end{lem}
\textbf{Proof.} Let $P$ be an essential minimal prime ideal, of
$C(X)$. Then by \cite[Corollary, $3.3$]{4.}, there is $p\in \beta
X\setminus I(X)$ such that $O^{p}\subseteq P$ so $O^{\beta
X\setminus I(X)}$ is contained in the intersection of essential
minimal prime ideals. Now let $f$ be an element of the
intersection of essential minimal prime ideals and $p\in \beta
X\setminus I(X)$. Then by \cite[Theorem, 3.1]{4.}, $O^{p}$ is an
essential ideal, which is the intersection of some essential
minimal prime ideals, therefore $f\in O^{p}$. Hence the
intersection of essential minimal prime ideals is contained in
$O^{\beta X\setminus I(X)}$. $\hfill\square$
%----------------------------------------------------------------------------
\\ \indent
An ideal $I$ of $R$ is called a pure ideal, if for any $f\in I$,
there is a $g\in I$ such that $f=fg$ (see \cite{1.}). The set of
all points in a topological space $X$ which have compact
neighborhoods is denoted by $X_{L}$. It is easily seen that
$X_{L}=coz(C_{K}(X))=\bigcup_{f\in C_{K}(X)}coz(f)$. Since $\beta
X\setminus X\subseteq \beta X\setminus I(X)$, we have,
$C_{F}(X)\subseteq O^{\beta X\setminus I(X)}\subseteq C_{K}(X),$
where $C_{K}(X)=O^{\beta X\setminus X},$  see \cite[7. F]{7.}. In
the following theorem we show that the intersection of essential
minimal prime ideals in $C(X)$ is equal to the socle of $C(X)$.
However; it need not be equal to the $C_{K}(X)$.
%-------------------------------------------------------------------
\begin{thm}\label{1.11}
\begin{itemize}
\item[\rm(i)] The intersection of all essential minimal prime
ideals of $C(X)$ is equal to the socle of C(X) (i.e.,
$C_{F}(X))$.\item[\rm(ii)] Every intersection of essential minimal
prime ideals of $C(X)$ is an essential ideal \ifif the set of
isolated points of $X$ is dense in $X$.\item[\rm(iii)] The
intersection of all essential minimal prime ideals of $C(X)$ is
equal to the $C_{K}(X)$ \ifif $I(X)=X_{L}=\bigcup_{f\in
C_{K}(X)}\overline{X\setminus Z(f)}$, i.e., $I(X)=X_{L}$ and
$C_{K}(X)$ is a pure ideal.
\end{itemize}
\end{thm}
\textbf{Proof.} (i) By Lemma \ref{1.10}, the intersection of
essential minimal prime ideals is $O^{\beta X\setminus I(X)}$.
Hence $C_{F}(X)\sub O^{\beta X\setminus I(X)}$. Now let $f\in
O^{\beta X\setminus I(X)}$. Then $\beta X\setminus int_{\beta
X}cl_{\beta X}Z(f)\sub I(X)$. By \cite[6.9 d]{7.} any isolated
point of $X$ is isolated in $\beta X$, so $\beta X\setminus
int_{\beta X}cl_{\beta X}Z(f)$ is a compact subset of $\beta X$
consisting of some isolated points. Therefore $\beta X\setminus
int_{\beta X}cl_{\beta X}Z(f)$ is finite, this implies that
$X\setminus Z(f)$ is finite. Thus $f\in C_{F}(X)$.\\ \indent (ii)
By (i), this is \cite[Corollary 3.3]{3.}.  \\\indent (iii) Let
$C_{K}(X)=O^{\beta X\setminus I(X)}$. It is easily seen that
$I(X)\subseteq X_{L}$. Now let $x\in X_{L}$, then there exists a
compact subset $U$ in $X$ such that $x\in intU$, i.e., $x\notin
X\setminus intU$. By completely regularity of $X$  there is $f\in
C(X)$ such that $x\in X\setminus Z(f)\subseteq U\subseteq
cl_{X}U$, hence $x\in X\setminus Z(f)$, where $f\in C_{K}(X)$. By
hypothesis, $X\setminus I(X)\subseteq Z(f)$ so $x\in I(X)$.
Therefore $I(X)=X_{L}$, hence $C_{K}(X)=O^{\beta X\setminus
I(X)}=O^{\beta X\setminus coz(C_{K}(X))}$. By [1, Theorem, $3.2$],
$X_{L}=coz(C_{K}(X))=\bigcup_{f\in C_{K}(X)}\overline{X\setminus
Z(f)}$.\\\indent Conversely, we have $I(X)=X_{L}=\bigcup_{f\in
C_{K}(X)}\overline{X\setminus Z(f)}$.  By \cite[Theorem, 3.2]{1.},
$C_{K}(X)=O^{\beta X\setminus coz(C_{K}(X))}=O^{\beta X\setminus
X_{L}}=O^{\beta X\setminus I(X)}$. $\hfill\square$
%----------------------------------------------------------
\\\indent Recall that a completely regular space $X$ is said to be a
pseudo-discrete space if every compact subset of $X$ has finite
interior. Clearly the class of pseudo-discrete spaces contains the
class of P-spaces (see \cite{3.}).
\begin{cor}
A topological space $X$ is pseudo-discrete  \ifif $I(X)=X_{L}$ and
$C_{K}(X)$ is a pure ideal.
\end{cor}
\textbf{Proof.} This is a consequence of Theorem 4.5 in \cite{3.}
and Theorem \ref{1.11}. $\hfill\square$
%-----------------------------------------------------------------------
\\\indent
By using the above theorem, we give examples of topological spaces
$X$ for which $C_{K}(X)$ is equal to the intersection of essential
minimal prime ideals (i.e., $X$ is a pseudo-discrete space).
%---------------------------------------------------------------------------------------
\begin{exam}
(i) If $X$ is a  locally compact space and $C_{K}(X)=O^{\beta
X\setminus I(X)}$, then $X$ is a discrete space. For if $X$ is a
locally compact, then $C_{K}(X)$ is a pure ideal and $X_{L}=X$ so
if $C_{K}(X)=O^{\beta X\setminus I(X)}$, then $I(X)=X_{L}=X$,
i.e., $X$ is a discrete space.
\\\indent (ii) Let $X=Q$ (i.e., the set of rational numbers), with all points
having their usual neighborhood  except for $x=0$ is isolated
point. Then $X_{L}=I(X)=\{0\}$ and $C_{K}(X)$=$\{f\in C(X)$: $f=0$
except for $x=0\}$ is a pure ideal so $C_{K}(X)=O^{\beta
X\setminus I(X)}$.
\end{exam}
%------------------------------------------------------------
In the following, we see an example of a space $X$ for which
$C_{K}(X)$ is not  equal to the intersection of essential minimal
prime ideals.
%--------------------------------------------------------------------------------------
\begin{exam}
Let $X=[-1, 1]$ with a topology in which $x=0$ has the usual
neighborhood topology and all  other points are isolated. Then
$X_{L}=X\setminus \{0\}=I(X)$ but $C_{K}(X)$ is not equal to
$O^{\beta X\setminus I(X)}$, for $C_{K}(X)$ is not a pure ideal,
see \cite[Example 3.3]{1.}.
\end{exam}

By Theorem \ref{1.8} and Theorem \ref{1.10}, we have the following
corollary.
%-----------------------------------------------------------------------
\begin{cor}
The set of isolated points of $X$ is dense in $X$ \ifif the set of
isolated points of $\mathcal{M}(C(X))$ is dense in
$\mathcal{M}(C(X))$.
\end{cor}
%-------------------------------------------------------------------
By \cite[Theorem, $2.4$]{3.}, and  Proposition \ref{1.9}, we have
the following corollary.
%----------------------------------------------------
\begin{cor}
Let $X$ be a compact space. Every first category subset of $X$ is
nowhere dense in $X$ \ifif every first category subset of
$\mathcal{M}(C(X))$ is nowhere dense in $\mathcal{M}(C(X))$.
\end{cor}
%--------------------------------------------------
\begin{ques}
\begin{itemize}
\item[\rm(i)] For $R=C(X)$,  determine  $X$ for which
$R_{K}(\mathcal{M})=C_{K}(X)$? Note that in case $X$ and
$\mathcal{M}$ are compact or nowhere compact, then
$R_{K}(\mathcal{M})=C_{K}(X)$.
 \item[\rm(ii)] When is the intersection of $sd$-ideals in a reduced ring $R$ a $sd$-ideal?
 \item[\rm(iii)] When is the intersection of $sd$-ideals in $C(X)$ a $sd$-ideal?
 \item[\rm(iv)] When is $R_{K}(\mathcal{M})$  equal to the intersection of all strongly dense $fd$-ideals?
 \end{itemize}
 \end{ques}

\begin{center}{\textbf{Acknowledgements}}
\end{center}
The author wish to thank Dr. E. Momtahen for his encouragement and
discussion on this paper. The author also would like to thank the
well-informed referee of this article for the detailed report,
corrections, and several constructive suggestions for improvement.
%------------------------------------------------------------------------------------%

%-----------------------------------------------------------------------------
%-----------------------------------------------------------------------------

\end{document}